\documentclass{IEEEtran4PSCC}
\ifCLASSINFOpdf
   \usepackage[pdftex]{graphicx}
\else
   \usepackage[dvips]{graphicx}
\fi
\usepackage[cmex10]{amsmath}
\usepackage{amssymb} 
\usepackage{xcolor}
\usepackage{textcomp}
\usepackage{lipsum}
\newcommand{\bm}{\boldsymbol}
\usepackage{booktabs}
\usepackage[subtle,tracking=normal]{savetrees}
\usepackage{microtype} 

\allowdisplaybreaks

\usepackage[ruled]{algorithm2e}
\usepackage{microtype} 
\usepackage[subtle,tracking=normal]{savetrees}
\usepackage{url}
\usepackage{hyperref} 
\usepackage{xcolor}

\hypersetup{
    colorlinks=true, 
    linkcolor=blue, 
    filecolor=black, 
    urlcolor=blue, 
    citecolor=blue, 
    linkbordercolor={1 1 1}, 
    pdfborder={0 0 0} 
}
\makeatletter
\let\old@ps@headings\ps@headings
\let\old@ps@IEEEtitlepagestyle\ps@IEEEtitlepagestyle
\def\psccfooter#1{%
    \def\ps@headings{%
        \old@ps@headings%
        \def\@oddfoot{\strut\hfill#1\hfill\strut}%
        \def\@evenfoot{\strut\hfill#1\hfill\strut}%
    }%
    \def\ps@IEEEtitlepagestyle{%
        \old@ps@IEEEtitlepagestyle%
        \def\@oddfoot{\strut\hfill#1\hfill\strut}%
        \def\@evenfoot{\strut\hfill#1\hfill\strut}%
    }%
    \ps@headings%
}
\makeatother

\psccfooter{%
        \parbox{\textwidth}{\hrulefill \\ \small{24th Power Systems Computation Conference} \hfill \begin{minipage}{0.2\textwidth}\centering \vspace*{4pt} \includegraphics[scale=0.06]{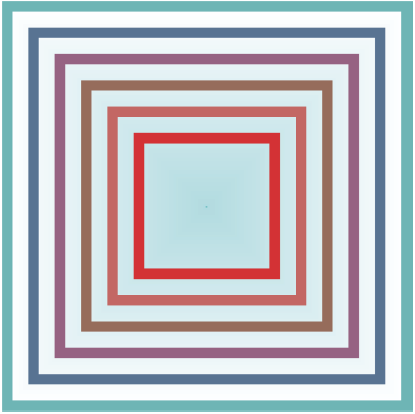}\\\small{PSCC 2026} \end{minipage} \hfill \small{Limassol, Cyprus --- June 8-12, 2026}}%
}

\begin{document}

\title{An Iterative Problem-Driven Scenario Reduction\\
Framework for Stochastic Optimization with \\Conditional Value-at-Risk}

\author{
  Yingrui Zhuang\textsuperscript{1*}, Lin Cheng\textsuperscript{1}, Ning Qi\textsuperscript{2*}, 
  Mads R. Almassalkhi\textsuperscript{3}, Feng Liu\textsuperscript{1}\\
    \textsuperscript{1}\textit{Department of Electrical Engineering, Tsinghua University, Beijing, China}\\
    \textsuperscript{2}\textit{Department of Earth and Environmental Engineering, Columbia University, New York, USA}\\
    \textsuperscript{3}\textit{Department of Electrical and Biomedical Engineering, University of Vermont, Burlington VT, USA}\\
    \textsuperscript{*}\textit{Email: zyr21@mails.tsinghua.edu.cn, nq2176@columbia.edu}
}


\maketitle
\begin{abstract}
Scenario reduction (SR) alleviates the computational complexity of scenario-based stochastic optimization 
with conditional value-at-risk (SBSO-CVaR)
by identifying representative scenarios to depict the underlying uncertainty and tail risks. 
Existing distribution-driven SR methods emphasize statistical similarity but often exclude extreme scenarios, 
leading to weak tail-risk awareness and insufficient problem-specific representativeness.
Instead, 
this paper 
proposes an iterative problem-driven scenario reduction framework.
Specifically, 
we integrate the SBSO-CVaR problem structure into SR process and project the original scenario set from distribution space onto the problem space.
Subsequently, 
to minimize the SR optimality gap with acceptable computation complexity,
we propose a tractable iterative problem-driven scenario reduction (IPDSR) method that
selects representative scenarios that best approximate the optimality distribution of the original scenario set while preserving tail risks.
Furthermore, 
the iteration process is rendered as a mixed-integer program to enable scenario partitioning and representative scenarios selection. 
And ex-post problem-driven evaluation indices are proposed to evaluate the SR performance.
Numerical experiments show IPDSR significantly outperforms existing SR methods by achieving an optimality gap of less than 1\% within an acceptable computation time.
\end{abstract}


\thanksto{\noindent This work is supported in part by National Natural Science Foundation of China (No. 52037006) and China Postdoctoral Science Foundation special funded project (No. 2023TQ0169).}

\section{Introduction}\label{sec:intro}

The increasing integration of renewable energy sources (RES) and emerging flexible loads introduces substantial uncertainties into power systems, 
significantly increasing operational risks~\cite{Uncertaintyreview}. 
Effective decision-making under uncertainty necessitates risk-averse optimization methods capable of accurately characterizing and quantifying these risks. 
Among risk-averse optimization frameworks such as robust optimization~\cite{zhang2021robust}, chance-constrained optimization~\cite{qi2023chance}, 
and distributionally robust optimization~\cite{wang2023two}, stochastic optimization with conditional value-at-risk (CVaR)~\cite{cvar1,zhu2024real} has been extensively applied in power systems, 
owing to convexity, coherency, and explicitly captures tail risks. Accurate CVaR computation requires precise modeling of uncertainty distributions, particularly in the tail regions. 
When exact distributional information is unavailable, discrete scenario approximations are commonly adopted. 
In this context, the scenario-based stochastic optimization with conditional value-at-risk (SBSO-CVaR)~\cite{cvar2} offers an effective way to manage tail risks. 
However, SBSO-CVaR relies on a large number of scenarios to accurately capture 
low-probability and high-loss events,
which leads to substantial computational complexity,
particularly as the dimensionality of uncertainty increases.
To alleviate computational burdens, scenario reduction (SR) is employed to identify a smaller representative scenario set to replace the
original scenario set while maintaining an acceptably robust optimal solution~\cite{SRreview}. 
Most existing SR techniques primarily focus on the risk-neutral SBSO formulations, but lack consideration for risk-averse formulations (e.g., SBSO-CVaR). 
Directly applying these methods to SBSO-CVaR is challenging because CVaR places particular emphasis on tail risk. 
Specifically, SR for SBSO-CVaR presents two main challenges:
\textit{i)} how to identify salient and representative scenarios that enable substantial scenario reduction, and
\textit{ii)} how to ensure that the reduced scenario set retains the ability to accurately capture tail risk behavior
and approximate the optimal solution of the original problem.

Existing SR methods for SBSO problems can be broadly classified into two categories:
\textit{distribution-driven scenario reduction} (DDSR) and \textit{problem-driven scenario reduction} (PDSR).
DDSR methods aim to preserve the statistical properties of the original scenario set, operating under the assumption that a more faithful representation of the underlying distribution will lead to a more accurate approximation of the original optimal solution.
Certain distribution-driven similarity metrics (e.g., Euclidean distance~\cite{qi2023portfolio}, 
Wasserstein distance~\cite{Wasserstein})
are commonly adopted by DDSR methods to measure the similarity between scenarios.
Distribution-driven clustering methods (e.g., K-means~\cite{k-means1}, hierarchical clustering~\cite{HC}, Gaussian mixture model~\cite{GMM})
and adopted to perform scenario clustering.
That is, DDSR methods generally consider SR and SBSO as two distinct and decoupled processes,
where the SR process is solely driven by statistical characteristics of scenarios.
However, 
since the same scenario may yield different solutions under different problem formulations 
(e.g., economic dispatch and resilience-oriented dispatch),
statistical similarity does not necessarily imply a similarity in the aspect of decision making~\cite{teichgraeber2019clustering}.
Therefore, DDSR methods suffer from a critical oversight: an inability to consider the impacts of the reduced scenario set on the
optimality.
To address this limitation, our recent work~\cite{PDSR} proposes a PDSR framework for risk-neutral SBSO problems. 
The key idea of PDSR is to project the original scenario set from the distribution space onto the problem space by explicitly incorporating the optimization problem structure. 
SR is then performed directly in the problem space by identifying scenarios that significantly influence decision-making optimality. 
The effectiveness of PDSR has been demonstrated through its ability to closely approximate the original optimal solution while maintaining computational efficiency.
However,
the PDSR framework in~\cite{PDSR} primarily focuses on risk-neutral problem and has not yet been extended to risk-averse settings.

Compared with risk-neutral SBSO problems, SBSO-CVaR problems exhibit greater complexity due to the presence of the CVaR term in their objective, 
requiring more precise modeling of tail risks~\cite{CVaRcomplexity}. 
Despite this additional complexity, 
most existing studies addressing SBSO-CVaR still employ traditional DDSR methods~\cite{SR_CVaR, LISRCVaR}. 
Such methods discards rare yet influential tail scenarios since they typically prioritize scenarios near the distribution center
and use averaging to derive representative scenarios.
Moreover, DDSR methods fail to account for the problem-specific characteristics of SBSO-CVaR or the sensitivity of tail risk measures to individual scenarios during SR. 
Consequently, these methods fail to capture the weighted marginal contributions of scenarios to tail losses, 
resulting in a poor representation of the tail risk distribution. 
This may beget different  decisions from those of the original problem, 
and ultimately, inadequate robustness to tail risks.

Existing PDSR methods for SBSO-CVaR focus exclusively on the CVaR term, which address tail risks but neglects other salient features of the underlying  distribution. 
For example, Ref.~\cite{SR_CVaR1} iteratively evaluates all tail scenarios, but results in a large number of scenarios and overlooks the impact of SR on the optimal solution. Consequently, this method suffers from scalability and suboptimality issues. In Ref.~\cite{SR_CVaR2}, tail-region scenarios capture the CVaR value well, but all remaining scenarios are aggregated into a single scenario. 
Such aggregation requires knowledge of the optimal solution under the original scenario set, 
which is generally unavailable in practice. 

Previous work on PDSR focused on risk-neutral SBSO problems~\cite{PDSR}. Directly applying this approach to risk-aware SBSO-CVaR settings can lead to inaccurate treatment of tail risks. In this paper, we extend the derivation process from~\cite{PDSR} to a CVaR-based formulation that captures risk-awareness and address the computational challenges that arise. Building on this framework, we introduce an iterative problem-driven scenario reduction (IPDSR) method that achieves near-optimal solutions to SBSO-CVaR problems with greatly reduced computational effort.
Specifically, our contributions are as follows:
\begin{enumerate}
  \item \textit{IPDSR Framework:} We propose an iterative problem-driven scenario reduction (IPDSR) framework that provide a tractable extension of~\cite{PDSR} to incorporate SBSO-CVaR. 
  IPDSR incorporates the problem structure of SBSO-CVaR (especially the emphasis on tail risks) into the SR process.  By analytically minimizing the optimality gap (OG), IPDSR 
  iteratively selects the representative scenarios and 
  ensures the reduced scenario set accurately approximates the original optimality.
    
  \item \textit{Solution Methodology:} 
  We iteratively identify the most representative scenarios within the projected problem space. In each iteration, the optimization process is reformulated as a mixed-integer programming (MIP) model that jointly performs scenario partitioning and representative scenario selection.
  Additionally, to efficiently handle large-scale problems, we propose an objective aggregation technique that significantly reduces the computational complexity of the MIP.
\item \textit{Simulation Analysis:} We evaluate the effectiveness of the proposed IPDSR framework using a risk-averse offering problem within the context of a virtual power plant (VPP) problem.
Simulation results demonstrate that IPDSR outperforms state-of-the-art SR methods by effectively identifying salient scenarios and achieving an optimality gap of less than 1\%. 
\end{enumerate}

The remainder of this paper is organized as follows. 
Section~\ref{sec:IPDSR} introduces the proposed IPDSR framework. Section~\ref{sec:formulation} presents the problem formulation of a risk-averse SBSO-CVaR problem. 
Section~\ref{sec:casestudy} presents case studies to validate
the effectiveness of IPDSR. 
We conclude the paper in Section~\ref{sec:conclusion}.

\section{Iterative Problem-Driven \\Scenario Reduction Framework}\label{sec:IPDSR}


\subsection{Formulation of SBSO-CVaR}
SBSO-CVaR represents a rich set of risk-averse power system problems.
The general formulation of SBSO-CVaR problem built on an original scenario set $\bm{\xi}$ is
{\small \begin{equation}\label{eq:CVaR_original}
    \min_{z\in Z}F(z,\bm{\xi}) \!:=\!  \sum_{i=1}^{N}\!\gamma_if(z,\xi_i) \!+\! \lambda \left( v_{\bm{\xi}}^{\alpha}\! + \!\frac{1}{1\!-\!\alpha} \sum_{i=1}^{N}\!\gamma_i  [ f(z,\xi_i)\! -\!v_{\bm{\xi}}^{\alpha} ]_{+} \right) ,
\end{equation}
}
where $F(z,\bm{\xi})$ is the objective function of the SBSO-CVaR problem with decision variable $z \in Z$,
including the expected cost across $\bm{\xi}$ (the first term) and the CVaR term (the second term) with $\alpha \in (0,1)$ as the risk preference. 
The risk aversion parameter is $\lambda \geq 0$, 
where $\lambda=0$ corresponds to the risk-neutral case and larger $\lambda$ indicates higher risk aversion.
$v_{\bm{\xi}}^{\alpha}$ is the VaR value under confidence level $\alpha$.
The positive part of $x$ is denoted $[x]_+ := \max\{x, 0\}$.
The inclusion of the CVaR term emphasizes the conditional expected loss in the tail distribution, transforming a risk-neutral optimization problem into a risk-averse one, which is inherently different from our previous work~\cite{PDSR}.
The original set of $N$ scenarios is denoted as $\bm{\xi}=\{\xi_1,\hdots, \xi_N\}$.
The probability of scenario $\xi_i$ is $\gamma_i > 0$, which satisfies $\sum_{i=1}^{N}\gamma_i=1$.
$f(z, \xi_i)$ denotes the scenario-specific cost function associated with scenario $\xi_i$.
We denote the optimal solution of the original problem~\eqref{eq:CVaR_original} as $z^\ast_{\bm{\xi}}= \underset{z\in Z}{{\arg\min}}\ F(z,\bm{\xi})$.
Similar to~\cite{PDSR},
we also reasonably require the SBSO-CVaR to satisfy the assumption of \textit{relatively complete recourse}.

Scenario reduction is employed to reduce the computational complexity by clustering large $N$ into $K\ll N$ clusters that preserve the accuracy of approximating optimality of the original problem.
The SR process can be denoted as 
$\bm{C}(\bm{\xi},K)=\{\{C_1,\hdots,C_K\}:C_i\neq  \emptyset, \forall i ;C_i\cap C_j=\emptyset,\forall i\neq j;\cup_i C_i=I, I = \{1,2,\hdots,N\}\}$. 
That is, the original scenario set $\bm{\xi}$ is partitioned into $K$ clusters
and is reduced to a representative scenario set $\bm{\zeta} = \{\zeta_1,\hdots, \zeta_K\}$. 
Each scenario cluster $C_k$ is represented by the representative scenario $\zeta_k$ with corresponding weight $\omega_k=\sum_{i \in C_k} \gamma _i$, which satisfies $\sum_{k=1}^{K}\omega_k=1$.
In this paper, we concentrate on selecting $\bm{\zeta} \subseteq \bm{\xi}$, 
instead of generating synthetic scenarios. 
The reduced problem built on the reduced scenario set $\bm{\zeta}$ can be formulated as~\eqref{eq:CVaR_reduced}
and the corresponding optimal solution is denoted as $z^\ast_{\bm{\zeta}}= \underset{z\in Z}{{\arg\min}}\ F(z,\bm{\zeta})$.

{\footnotesize
\begin{equation}\label{eq:CVaR_reduced}
    {
    \min_{z\in Z} F(z,\bm{\zeta}) \!:= \!\sum_{k=1}^{K}\omega_k  \!f(z,\zeta_k) \!
    + \!\lambda \left( \!v_{\alpha}^{\bm{\zeta}} \!+\! \frac{1}{1\!-\!\alpha}\sum_{k=1}^{K} \!\omega_k \!\left[f(z,\zeta_k)\!-\!v_{\alpha}^{\bm{\zeta}}\right]_{+}\right).
    }
\end{equation}
}
CVaR emphasizes tail risk characterization, inducing non-smooth piecewise convex structures that markedly increase computational complexity. As a result, SR for SBSO-CVaR poses substantial challenges beyond the risk-neutral case.
In what follows, we extend the framework to accommodate the SBSO-CVaR setting.


\subsection{Optimality Gap of Scenario Reduction}
SR aims to minimize the optimality gap (OG) from using $K$ representatives ($\bm{\zeta}$) vs. $N$ scenarios ($\bm{\xi}$). 
Towards this purpose, the OG can be defined as
\begin{equation}\label{eq:optimality gap}
OG:=F(z^\ast_{\bm{\zeta}},\bm{\xi})-F(z^\ast_{\bm{\xi}},\bm{\xi}),
\end{equation}
where $F (z^\ast_{\bm{\zeta}},\bm{\xi})$ means solving~\eqref{eq:CVaR_original} with $z=z^\ast_{\bm{\zeta}}$.
Since $F(z,\bm{\xi})\ge F(z^\ast_{\bm{\xi}},\bm{\xi})$ for all $z\in Z$, we have $OG \ge 0$.
A smaller $OG$ indicates a more accurate problem optimality approximation of $\bm{\zeta}$ to $\bm{\xi}$. 

Next, we present two methods to minimize OG:

\textit{i)} PDSR: 
We extend~\cite{PDSR} for risk-neutral formulations to risk-averse formulations with CVaR,
and derive a theoretical upper bound of the OG leveraging relaxation techniques,
as detailed in Section~\ref{subsec:PDSR} and Appendix~\ref{sec:appendix_PDSR}.

\textit{ii)} IPDSR: we iteratively select the representative scenario set $\bm{\zeta}$ to minimize the OG, as detailed in Section~\ref{subsec:IPDSR}.

\subsection{Problem-Driven Scenario Reduction}\label{subsec:PDSR}
This method provides a risk-averse formulations with CVaR,
and further incorporates the problem structure of SBSO-CVaR 
to devise a relaxation strategy. 
We first  project the distributional space onto the problem space and define a problem-driven distance metric, 
upon which the SR process is reformulated as an MIP model aiming to minimize the relaxed theoretical upper bound on the OG.
The detailed derivations are presented in Appendix~\ref{sec:appendix_PDSR}.

However, the MIP formulation 
involves numerous binary variables and  constraints, 
making it computationally intractable for standard solvers. Therefore, 
we introduce an alternative approach to minimize the OG and greatly reduce the SR problem complexity.
Since $z^\ast_{\bm{\xi}}$ remains unknown in practice and $F(z^\ast_{\bm{\xi}},\bm{\xi})$ is constant, 
the OG can be minimized by minimizing $F(z^\ast_{\bm{\zeta}},\bm{\xi})$.
We next develop an IPDSR framework to 
identify a good $\bm{\zeta}$ and corresponding $z^\ast_{\bm{\zeta}}$
that minimize the OG.

\subsection{Iterative Problem-Driven Scenario Reduction}\label{subsec:IPDSR}
IPDSR contains two main steps:
initialization and iterative scenario partition.

\subsubsection{Initialization}
Firstly, $\bm{\zeta}^0$ and $z_{\bm{\zeta}^0}^\ast$ are initialized.
We denote $\kappa  $ as the iteration index, initialized as $\kappa =0$. 
We validate $z_{\bm{\zeta}^0}^\ast$ in $\bm{\xi}$ 
and generate validation value of $F(z_{\bm{\zeta}^0}^\ast,\bm{\xi})$
and the initial VaR value $v_{\bm{\xi}}^{\alpha,0}$.
Note that this initialization can be done by applying traditional DDSR methods like $K$-means.
And multiple initializations can be performed to obtain a better initial $\bm{\zeta}^0$.


\subsubsection{Iterative Scenario Reduction}\label{subsec:ISR}
At step $\kappa$, we have $\bm{\zeta}^{\kappa }$, $z_{\bm{\zeta}^{\kappa }}^\ast$, $F(z_{\bm{\zeta}^{\kappa}}^\ast,\bm{\xi})$ 
and $v_{\bm{\xi}}^{\alpha,\kappa}$. 
We seek to select a new reweighted representative scenario set $\bm{\zeta}^{\kappa +1}\subset \bm{\xi}$ with weights $\omega^{\kappa +1}$ 
that $F(z^\ast_{\bm{\zeta}^{\kappa }},\bm{\zeta}^{\kappa +1})$ can best approximate $F(z_{\bm{\zeta}^{\kappa}}^\ast,\bm{\xi})$,
which can be formulated as the following optimization problem:
\begin{align}
&\bm{\zeta}^{\kappa+1}  = \underset{\bm{\zeta} \subset \bm{\xi}}{\arg\min} 
   | F(z^\ast_{\bm{\zeta}^\kappa}, \bm{\zeta}) - F(z^\ast_{\bm{\zeta}^\kappa}, \bm{\xi}) | 
   = \underset{\bm{\zeta} \subset \bm{\xi}}{\arg\min}\notag \\
& 
   | \sum_{k=1}^{K} \sum_{i \in C_k} \gamma_i ( 
   f(z^\ast_{\bm{\zeta}^\kappa}, \zeta_k) 
   + \lambda ( v_{\bm{\zeta}}^{\alpha, \kappa} 
   + \frac{1}{1 - \alpha} \gamma_i 
   [ f(z^\ast_{\bm{\zeta}^\kappa}, \zeta_k) 
   - v_{\bm{\zeta}}^{\alpha, \kappa} ]_+ 
   ) )  \notag \\
& 
   - \sum_{k=1}^{K} \sum_{i \in C_k} \gamma_i ( 
   f(z^\ast_{\bm{\zeta}^\kappa}, \xi_i) 
   + \lambda ( v_{\bm{\xi}}^{\alpha, \kappa} 
   + \frac{1}{1 - \alpha} \gamma_i 
   [ f(z^\ast_{\bm{\zeta}^\kappa}, \xi_i) 
   - v_{\bm{\xi}}^{\alpha, \kappa} ]_+ 
   ) ) | \notag \\
&= \underset{\bm{\zeta} \subset \bm{\xi}}{\arg\min} 
   \left| \sum_{k=1}^{K} \sum_{i \in C_k} \gamma_i \Big( 
   f(z^\ast_{\bm{\zeta}^\kappa}, \xi_i) - f(z^\ast_{\bm{\zeta}}, \zeta_k) 
   + \lambda (v_{\bm{\xi}}^{\alpha, \kappa} - v_{\bm{\zeta}}^{\alpha, \kappa}) \right. \notag \\
& 
   \left. + \frac{\lambda}{1 - \alpha} [ f(z^\ast_{\bm{\zeta}}, \xi_i) 
   - v_{\bm{\xi}}^{\alpha, \kappa} ]_+ 
   - [ f(z^\ast_{\bm{\zeta}}, \zeta_k) 
   - v_{\bm{\zeta}}^{\alpha, \kappa} ]_+ 
   \Big) \right|. \label{eq:Iter_obj}
\end{align}

For convenience of notation, we omit the superscript $\kappa$ in the following discussion.
We sort the scenarios in non-decreasing manner with respect to the validated objective function $f(z^\ast_{\bm{\zeta}},\cdot)$, i.e.,
$f(z^\ast_{\bm{\zeta}},\xi_{(1)})\leq f(z^\ast_{\bm{\zeta}},\xi_{(2)})\leq \dots \leq f(z^\ast_{\bm{\zeta}},\xi_{(N)})$,
where $\xi_{(i)}$ is the ordered $i$-th scenario with corresponding probability $\gamma_{(i)}$.
Let the integer index $N^{\alpha}_{\bm{\xi}}$ be defined such that
$\sum_{i=1}^{N^{\alpha}_{\bm{\xi}}-1}\gamma_{(i)} < \alpha \leq \sum_{i=1}^{N^{\alpha}_{\bm{\xi}}}\gamma_{(i)}$.
Then we have $v_{\bm{\xi}}^{\alpha} = f(z,\xi_{(N^{\alpha}_{\bm{\xi}})})$,
and for $i > N^{\alpha}_{\bm{\xi}}$, we have $\left[ f(z, \xi_{(i)}) - v_{\bm{\xi}}^{\alpha} \right]_{+} = f(z, \xi_{(i)}) - v_{\bm{\xi}}^{\alpha}$,
otherwise 0.
Similar for $\bm{\zeta}$, we have $K^{\alpha}_{\bm{\zeta}}$ such that  
$v_{\bm{\zeta}}^{\alpha} = f(z^\ast_{\bm{\zeta}},\zeta_{(K^{\alpha}_{\bm{\zeta}})})$.

In this way, we project the original scenario set $\bm{\xi}$ from the distribution space
onto the problem space of SBSO-CVaR by explicitly incorporating the problem structure.
Note that the problem space here is a collapsed space since we have the validation
solution $z^\ast_{\bm{\zeta}}$.
The projection process can be denoted as 
$\boldsymbol{\xi}\rightarrow \boldsymbol{F}$,
where $\boldsymbol{F}:=\{F_{i}=f(z^\ast_{\bm{\zeta}},\xi_{(i)})\}$.
In this projection, we can directly quantify the impacts of uncertainty 
and systematically incorporate the inherent characteristics of the SBSO-CVaR problem into the SR process.

Then, we transform the optimization problem~\eqref{eq:Iter_obj} into a mixed-integer programming problem as~\eqref{eq:mip},
which conducts scenario partition and representative scenario selection in the problem space simultaneously.
\begin{subequations}\label{eq:mip}
    \begin{align}
        \min
        &\ |\sum_{j=1}^{N}(\sum_{i=1}^{N}v_{ij}\gamma_i(F_i-F_j+\frac{\lambda}{1-\alpha}([F_i-v_{\bm{\xi}}^{\alpha}]_+ - [F_j - v_{\bm{\zeta}}^{\alpha}]_+)))\notag
                \\&\ \ \ \ \ +\lambda(v_{\bm{\xi}}^{\alpha}-v_{\bm{\zeta}}^{\alpha})| \label{eq:mip_obj} \\
        \text{s.t. } 
        & v_{ij} \leq u_{j},\ v_{jj} = u_{j},\ n_j \leq u_j,\ h_j \leq u_j, \\
        & \sum\nolimits_{j=1}^{N}v_{ij} = 1,\ \sum\nolimits_{j=1}^{N}u_{j} = K, \label{eq:assignment} \\
        & \sum\nolimits_{j=1}^{N}n_j = 1, \ \sum\nolimits_{i=1}^{N}n_i F_i =v_{\bm{\zeta}}^{\alpha}, \label{eq:var1} \\
        & h_{i} = u_{i}(1-\sum\nolimits_{p=1}^{i}n_{p}),\label{eq:varleftright}\\
        & \sum\nolimits_{j=1}^{N}h_j\sum\nolimits_{i=1}^{N}v_{ij}\gamma_i \leq \alpha, \label{eq:varprob1} \\
        & \sum_{j=1}^{N}h_j\sum_{i=1}^{N}v_{ij}\gamma_i + (u_k - h_k)\sum_{i=1}^{N}v_{ik}\gamma_i \geq \alpha (u_k - h_k), \label{eq:varprob2} \\
        & v_{ij} \in \{0,1\} ,\ u_{j} \in \{0,1\} ,\ n_{j} \in \{0,1\} ,\ h_{j} \in \{0,1\}. \label{eq:constr_binary_cont}
    \end{align}
\end{subequations}
where 
$u_j$ indicates whether scenario $\xi_j$ is selected as a representative scenario of a cluster.
$v_{ij}$ determines whether scenario $\xi_i$ is included in the cluster represented by scenario $\xi_j$. 
$n_j$ indicates whether scenario $\xi_j$ is the scenario corresponding to the VaR value (i.e., whether $j=K^{\alpha}_{\bm{\zeta}}$),
$h_j$ indicates $\mathbf 1_{j\leq K^{\alpha}_{\bm{\zeta}}}$ and~\eqref{eq:varleftright} is only valid for sorted scenarios.
Constraint $v_{ij}\leq u_{j}$ ensures that $\xi_i$ can only be assigned to a cluster that has a designated representative, 
while $v_{jj} = u_j$ enforces that $\xi_i$ must be assigned to its respective cluster if it's a representative scenario.
Constraint $n_j \leq u_j$ ensures that the VaR scenario of $\bm{\zeta}$ can only be selected from the representative scenarios.
Constraint $\sum_{j=1}^{N}v_{ij}=1$ ensures that each scenario $\xi_j$ can only be assigned to one cluster,
while $\sum_{j=1}^{N}u_{j}=K$ guarantees that exactly $K$ clusters are formed. 
Constraint $\sum_{j=1}^{N}n_j=1$ ensures that there is only one VaR scenario in the reduced scenario set $\bm{\zeta}$.
Constraints~\eqref{eq:varprob1} and~\eqref{eq:varprob2} are used to ensure the confidence level of the VaR scenario in the reduced scenario set $\bm{\zeta}$.

Note that the MIP model is nonlinear due to the presence of absolute values and variable products. 
Auxiliary variables can be introduced to partially linearize the model to make it solvable by Gurobi.
By solving the MIP model in~\eqref{eq:mip}, we can obtain the optimal representative scenario set 
$\bm{\zeta}^{\kappa+1}$ and the associated weights $\omega_k^{\kappa+1} = \sum_{i=1}^{N}v_{ik}^{\kappa+1} \gamma_{(i)}$.
Accordingly, the corresponding solution $z_{\bm{\zeta}^{\kappa+1}}^\ast$ is computed.


\subsection{Objective Aggregation for MIP Simplification}
Notably, the MIP model in~\eqref{eq:mip} is computationally intensive with lots of binary variables, especially when $N$ is large.
To address this issue, we propose to precede the solution process~\eqref{eq:mip} with an objective aggregation step. 
Specifically, $\boldsymbol{F}=\{F_{i}\}_{i=1}^{N}$ is aggregated into $\bar{\boldsymbol{F}}:=\{\bar{F}_{j}\}_{j=1}^{N'}$,
and the corresponding probability weights $\{\gamma_i\}_{i=1}^{N}$ are aggregated into $\{\bar{\gamma}_{j}\}_{j=1}^{N'}$,
where $N' \ll N$ is the number of aggregated values.
\begin{align}
\bar{F}_j = \frac{\sum_{i \in \mathcal{I}_j} \gamma_i F_i}{\sum_{i \in \mathcal{I}_j} \gamma_i}, 
\quad
\bar{\gamma}_j = \sum_{i \in \mathcal{I}_j} \gamma_i.
\end{align}
where $\mathcal{I}_j$ is the index set of aggregated cluster $j$.
These aggregated values are then used in the MIP model~\eqref{eq:mip} in place of the original ones, 
effectively reducing model size while preserving the key characteristics of the scenario distribution. 
Then, the results of the MIP model are re-mapped to the original scenario set $\bm{\xi}$.


\subsection{Algorithm of IPDSR}
The overall IPDSR algorithm is summarized in Algorithm~\ref{algo:IPDSR}.
\begin{algorithm}[htbp]\label{algo:IPDSR}
\caption{IPDSR}
\SetAlgoLined
\SetEndCharOfAlgoLine{}
\KwIn{Scenario set $\bm{\xi}$ of $N$ scenarios and weights $\bm{\gamma}$.}
\KwOut{Scenario set $\bm{\zeta}$ of $K$ scenarios and weights $\bm{\omega}$.}
\SetKwBlock{StepOne}{Step 1 - Initialization}{}
\SetKwBlock{StepTwo}{Step 2 - Iterative Scenario Reduction}{}

\StepOne{
    Initialize  $\bm{\zeta}^0  $ with $z_{\bm{\zeta}^0 }^\ast$ and $\kappa = 0$.\;
}
\StepTwo{
    \While{$\kappa  < \text{max\_iter}$} {
    Project the scenario set $\bm{\xi}$ in the distribution space onto the problem space
    by computing $\bm{F}$.\;
    Perform objective aggregation on $\bm{F}$ and \newline generate $\bm{\bar{F}}$ for MIP simplification.\;
    Construct and solve the MIP model in~\eqref{eq:mip} on $\bm{\bar{F}}$ \newline to obtain $\bm{\zeta}^{\kappa +1}$ with weights $\bm{\omega}^{\kappa +1}$.\;
    Update $z_{\bm{\zeta}^{\kappa +1}}^\ast = \underset{z\in Z}{{\arg\min}}\ F(z,\bm{\zeta}^{\kappa +1})$.\;
    Compute the validation value of $F(z_{\bm{\zeta}^{\kappa +1}}^\ast,\bm{\xi})$.\;
    Set $\kappa  = \kappa  + 1$.\;
    }
    Select $\bm{\zeta}^{\kappa +1}$ and corresponding $\bm{\omega}^{\kappa +1}$
    with smallest \newline $F(z_{\bm{\zeta}^{\kappa +1}}^\ast,\bm{\xi})$  
    as the final output.\;
}

\end{algorithm}

\subsection{Problem-Driven Evaluation Indices}\label{eval metric}
In this part, we introduce three \textit{Ex-post} indices focus on the impacts of SR on the optimality of SBSO-CVaR after solving the reduced problem.

\subsubsection{Validation Distribution Similarity}\label{DS}
The distribution similarity between $\{f(z^\ast_{\bm{\zeta}},\xi_{(i)})\}_{i=1}^N$ 
and $\{f(z^\ast_{\bm{\xi}},\xi_{(i)})\}_{i=1}^N$ can be evaluated by the Wasserstein distance (WD):
\begin{equation}\label{eval:Wasserstein}
    WD(\bm{\zeta},\bm{\xi}) = \sum\nolimits_{i=1}^{N} \gamma_i \left|f(z^\ast_{\bm{\zeta}},\xi_{(i)}) - f(z^\ast_{\bm{\xi}},\xi_{(i)})\right|,
\end{equation}
where $\xi_{(i)}$ is the $i$-th scenario in the ordered set.
A smaller WD implies a closer match between the validation outcomes, 
indicating that the reduced scenario set yields a solution with stronger applicability to the original scenarios.

\subsubsection{Optimality gap}\label{OG}
After solving~\eqref{eq:CVaR_original} and~\eqref{eq:CVaR_reduced}, 
we apply $z^\ast_{\bm{\xi}}$ and $z^\ast_{\bm{\zeta}}$ to $\bm{\xi}$, 
and calculate the percentage value of OG as
\begin{equation}\label{OG in eval index}
    OG_{\bm{\zeta}} (\%) :=\frac{F(z^\ast_{\bm{\zeta}},\bm{\xi})-F(z^\ast_{\bm{\xi}},\bm{\xi})}{F(z^\ast_{\bm{\xi}},\bm{\xi})}.
\end{equation}

This metric indicates the percentage derivation of the approximated optimality from the optimality of the original problem, and is desired to be close to zero.

\subsubsection{Representative scenario effectiveness}\label{effective scenario}
For the SBSO with representative scenarios, we introduce Scenario Effectiveness ($SE_{\zeta_k}(\%)$) to quantify the importance of scenario $\zeta_k$. It is measured by the change in the percentaged $OG$ when $\zeta_k$ is removed from $\bm{\zeta}$:
\begin{equation}
    SE_{\zeta_k}(\%) :=OG_{\bm{\zeta}_{-k}}(\%)-OG_{\bm{\zeta}}(\%),
\end{equation}
where $\bm{\zeta}_{-k}=\bm{\zeta}\setminus{\zeta_k}$.
A larger $SE_{\zeta_k}(\%)$ indicates greater influence of $\zeta_k$ on the reduced problem outcomes.
More details of this metric can be found in~\cite{PDSR}.

\section{Day-Ahead Risk-Averse Offering Problem Formulation of a Virtual Power Plant}\label{sec:formulation}

We consider a day-ahead risk-averse offering problem for a virtual power plant (VPP) trading with the upper grid. 
The VPP aggregates renewable energy sources (RES) and energy storage (ES) facilities, 
managing uncertainties from net load demand and market price fluctuations,
which are described using the representative scenario set. 
The VPP makes risk-averse decisions on day-ahead power trading while considering intra-day adjustments for derivations between the scheduled day-ahead power trading and the actual intraday power demand. 

The VPP day-ahead risk-averse offering problem is formulated in~\eqref{problem}. 
\begin{subequations}\label{problem}
  \begin{align}
    & \min\limits_{\Xi} \sum\nolimits_{i=1}^{N}\gamma_iC_i + \lambda (v_{\alpha} + \frac{1}{1-\alpha}\sum\nolimits_{i=1}^{N}\gamma_i  [ C_i - v_{\alpha} ]_{+} ) \label{obj_all}\\
    & C_{i} = \Delta t\sum\nolimits_{t=1}^{T} (\pi_{i,t}^{\mathrm{T}}P_{t}^{\mathrm{T}}+ 1.3\pi_{i,t}^{\mathrm{T}}P_{i,t}^{\mathrm{T+}}-0.7\pi_{i,t}^{\mathrm{T}}P_{i,t}^{\mathrm{T-}})\label{obj}\\
    & SoC_{i,t+1} = SoC_{i,t} +  t(P_{i,t}^{\mathrm{E},\mathrm{c}}\eta^{\mathrm{c}}-P_{i,t}^{\mathrm{E},\mathrm{d}}/\eta^{\mathrm{d}})/E \label{SoC0_N}\\
    & \sum\nolimits_{t=1}^{T}(P_{i,t}^{\mathrm{E},\mathrm{c}}\eta^{\mathrm{c}}-P_{i,t}^{\mathrm{E},\mathrm{d}}/\eta^{\mathrm{d}}) \Delta t = 0 \label{ES_ch_dis}\\
    & 0 \leq P_{i,t}^{\mathrm{E},\mathrm{c}} \leq (1-D_{i,t}^{\mathrm{E}})\overline{P^{\mathrm{E}}},\ 0 \leq P_{i,t}^{\mathrm{E},\mathrm{d}} \leq D_{i,t}^{\mathrm{E}}\overline{P^{\mathrm{E}}} \label{SoC_change} \\    
    & \underline{SoC} \leq SoC_{i,t} \leq \overline{SoC} \label{SoC_range}\\
    & P_{t}^{\mathrm{T}}\leq \overline{P}_{t}^{\mathrm{T}},\ -\overline{P}_{t}^{\mathrm{T}}\leq P_{t}^{\mathrm{T}}+P_{i,t}^{\mathrm{T+}}-P_{i,t}^{\mathrm{T-}}\leq \overline{P}_{t}^{\mathrm{T}}\label{PEX_range}\\
    & 0\leq P_{i,t}^{\mathrm{R,u}}\leq P_{i,t}^{\mathrm{R}}\label{PR_range}\\
    & -P_{i,t}^{\mathrm{R,u}}+P_{i,t}^{\mathrm{L}}+P_{i,t}^{\mathrm{E,c}}-P_{i,t}^{\mathrm{E,d}}=P_{t}^{\mathrm{T}}+P_{i,t}^{\mathrm{T+}}-P_{i,t}^{\mathrm{T-}}\label{P_balance}\\
    & 0 \leq P_{i,t}^{\mathrm{T+}} \leq D_{i,t}^{\mathrm{T}}\overline{P}^{\mathrm{T+}},\ 0 \leq P_{i,t}^{\mathrm{T-}} \leq (1-D_{i,t}^{\mathrm{T}})\overline{P}^{\mathrm{T-}}\label{PEX_state}\\
    & D_{i,t}^{\mathrm{E}},\ D_{i,t}^{\mathrm{T}} \in \{0,1\}\text{,}
  \end{align}
\end{subequations}
where $C_i$ is the cost of scenario $i$,
$\pi_{i,t}^{\mathrm{T}}$ is the day-ahead market price of scenario $i$ at time step $t$,
$P_{i,t}^{\mathrm{T}}$ is the day-ahead power trading,
$P_{i,t}^{\mathrm{T+}}/P_{i,t}^{\mathrm{T-}}$ are the up/down intraday balancing market power,
$P_{i,t}^{\mathrm{E},\mathrm{c}}/P_{i,t}^{\mathrm{E},\mathrm{d}}$ are the charging/discharging power of the ES,
$E/SoC_{i,t}$ are the energy capacity/state of charge of the ES,
$\eta^{\mathrm{c}}/\eta^{\mathrm{d}}$ are the charging/discharging efficiency of the ES,
$D_{i,t}^{\mathrm{E}}$ and $D_{i,t}^{\mathrm{T}}$ are the binary variables indicating the ES charging state and the intraday balancing state, respectively.
$\overline{P^{\mathrm{E}}}$ is the power rating of the ES,
$\overline{P}_{t}^{\mathrm{T}}$ is the maximum power exchange with the upper network.
$\underline{SoC}/\overline{SoC}$ are the minimum and maximum state of charge of the ES, respectively.
$\Xi$ is the set of decision variables.
Eq.~\eqref{obj} is the scenario-specific economic cost of the VPP,
including the cost of day-ahead trading and intraday power balancing.
Constraint~\eqref{ES_ch_dis} limits the charging and discharging power and state of the ES.
Constraint~\eqref{SoC_range} bounds the range of state of charge (SoC). 
Constraint~\eqref{SoC0_N} guarantees that the capacity at the last time period is equal to the initial capacity.
Constraint~\eqref{PEX_range} limits the day-ahead power trading and intraday balancing power.
Constraint~\eqref{PR_range} limits the renewable power curtailment.
Constraint~\eqref{P_balance} is the power balance equation of the VPP.
Constraint~\eqref{PEX_state} ensures that the VPP cannot simultaneously buy and sell power in the intraday balancing market.

\section{Numerical Case Study}\label{sec:casestudy}

\subsection{Set up}
In the considered VPP operation problem,
each individual scenario $\xi\in \mathbb{R} ^{2\times T},~\xi \in \boldsymbol{\xi}$ is a multi-variable high-dimensional vector 
characterizing 2~sources of uncertainty (net load demand and market price).
The time step is set as $\Delta t = 15\min$ with $T=96$.
The capacity of the wind power plant is set as $1000\mathrm{kW}$,
the capacity of the ES is set as $200\mathrm{kW}$ with 2 hours of charging/discharging duration.
The initial $SoC$ of ES is $0.5$ and $\eta^{\mathrm{c}}=\eta^{\mathrm{d}}=0.95$.
We set $\lambda = 0.5$ and $\alpha = 0.95$.
The problem built on $\bm{\xi}$ is used as the Benchmark. 
The optimization is coded in Python with the Gurobipy interface and solved by Gurobi 12.0 solver.
The programming environment is Intel Core i9-13900HX @ 2.30GHz with RAM 16 GB.
         
\subsection{Performance of IPDSR}

\begin{figure}[htbp] 
    \setlength{\abovecaptionskip}{-0.1cm}  
    \setlength{\belowcaptionskip}{-0.1cm}   
    \centerline{\includegraphics[width=0.95\columnwidth]{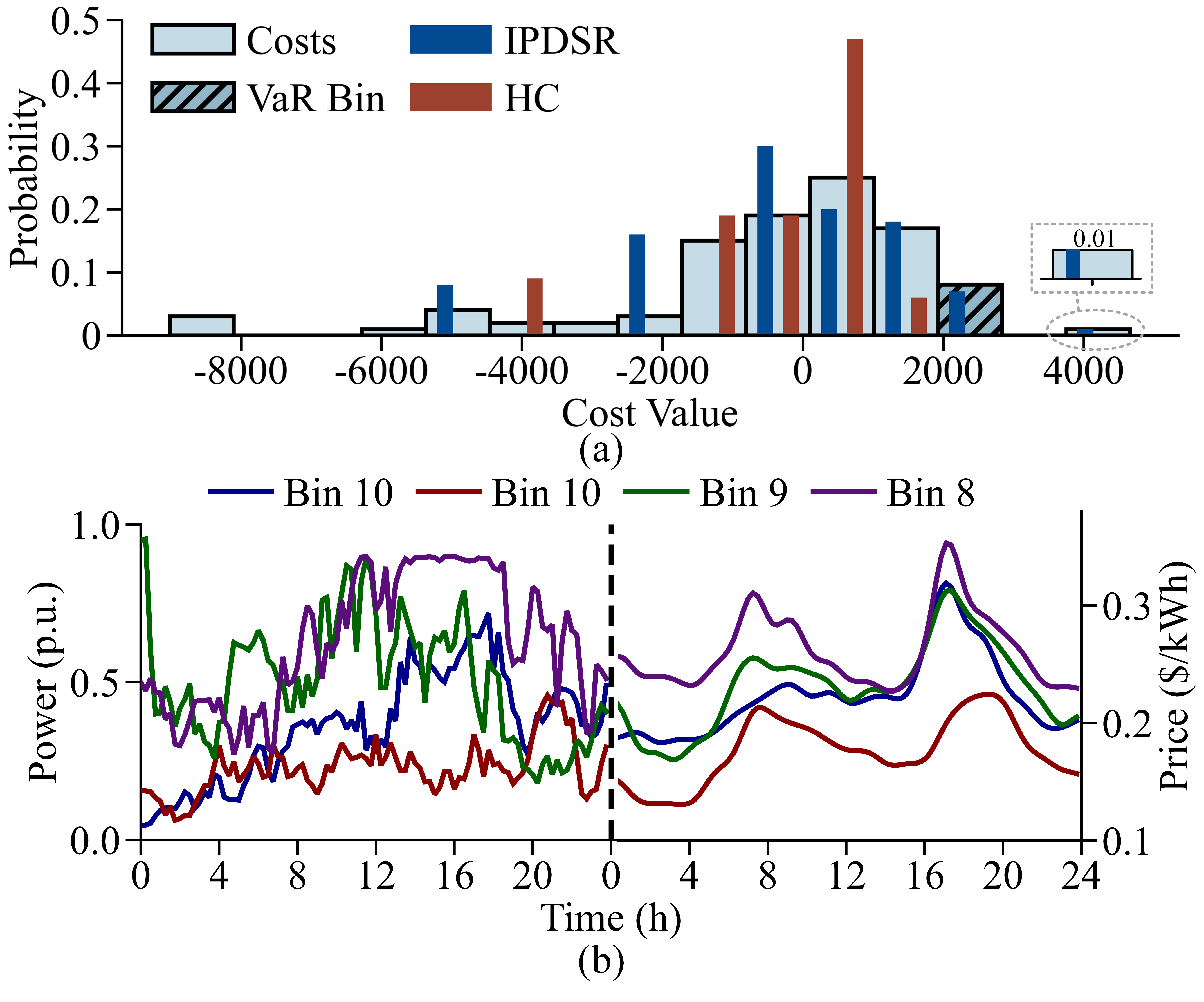}}
    \caption{Illustration example of non-linear mapping between distribution space and problem space:
    (a) scenarios in the problem space;
    (b) four example scenarios in the distribution space: left – net load scenarios; right – corresponding price scenarios.}
    \label{PDSR_DDSR} 
\end{figure}

We first visualize the nonlinear mapping between the distribution space and the problem space
to highlight the effectiveness of IPDSR and the limitations of classical DDSR methods 
(e.g., hierarchical clustering, HC). 
For clarity, we consider a small-scale example with $N = 100$ and $K = 10$.
We solve the original optimization problem~\eqref{eq:CVaR_original} over $\bm{\xi}$ to obtain $z^\ast_{\bm{\xi}}$ 
and then compute the scenario-specific objective values $\{f(z^\ast_{\bm{\xi}}, \xi_i)\}_{i=1}^N$. 
The histogram of $\{f(z^\ast_{\bm{\xi}}, \xi_i)\}_{i=1}^N$ is shown as the light-colored bars in Fig.~\ref{PDSR_DDSR}(a),
with the VaR threshold indicated by a shaded bar. 
Fig.~\ref{PDSR_DDSR}(a) illustrates the distribution of all scenarios in the problem space.
Since the decision $z^\ast_{\bm{\xi}}$ is given, 
the problem space collapses to the one-dimensional space of values. 
Notably, there exists an extreme scenario in the right tail of the distribution,
which is critical for SBSO-CVaR problems.

We then utilize the IPDSR framework to obtain the reduced scenario set $\bm{\zeta}$ and associated weights $\bm{\omega}$. 
Values of $\{f(z^\ast_{\bm{\xi}}, \zeta_k)\}_{k=1}^K$ are then mapped onto the same histogram as narrow blue bars, 
where scenarios falling into the same bin have their probabilities aggregated.
For comparison with DDSR methods, 
we also mark the representative scenarios obtained by the HC method as narrow orange bars using the same procedure.
Fig.~\ref{PDSR_DDSR}(a) shows that scenarios obtained by IPDSR are more widely spread across the objective value distribution 
and effectively capture the objective value distribution features, especially in the right tail.
In contrast, HC-selected scenarios are concentrated in the central region, 
failing to represent right-tail outcomes. 
This contrast stems from the methodological difference: 
DDSR performs clustering in the distribution space, 
while IPDSR clusters in the risk-aware problem space. 
Due to the nonlinear relationship between the distribution and problem spaces, 
similarity in distribution features does not guarantee effective decision-making in the problem space.

This point is further illustrated in Fig.~\ref{PDSR_DDSR}(b), where we present four scenarios samples. 
Although the blue and red scenarios exhibit significant differences in the distribution space, 
they fall into the same bin in the problem space. 
Conversely, the green and purple scenarios are distributionally similar but belong to different bins in the problem space. 
This clearly indicates that distribution similarity does not necessarily translate into decision-relevant similarity.

Moreover, 
the optimality gap of IPDSR is only $0.50\%$, 
whereas that of HC is $6.83\%$, 
highlighting the superior performance of IPDSR. 
This is largely because IPDSR accounts for SBSO-CVaR problem structure (particularly tail risk) in the SR process,
leading to more robust and adaptable representative scenarios and solutions.
Overall, these results highlight the advantage of IPDSR in preserving tail-risk problem structures 
and improving the representativeness of reduced scenario sets in the problem space.

\begin{figure}[htbp] 
    \setlength{\abovecaptionskip}{-0.1cm}  
    \setlength{\belowcaptionskip}{-0.1cm}   
    \centerline{\includegraphics[width=0.95\columnwidth]{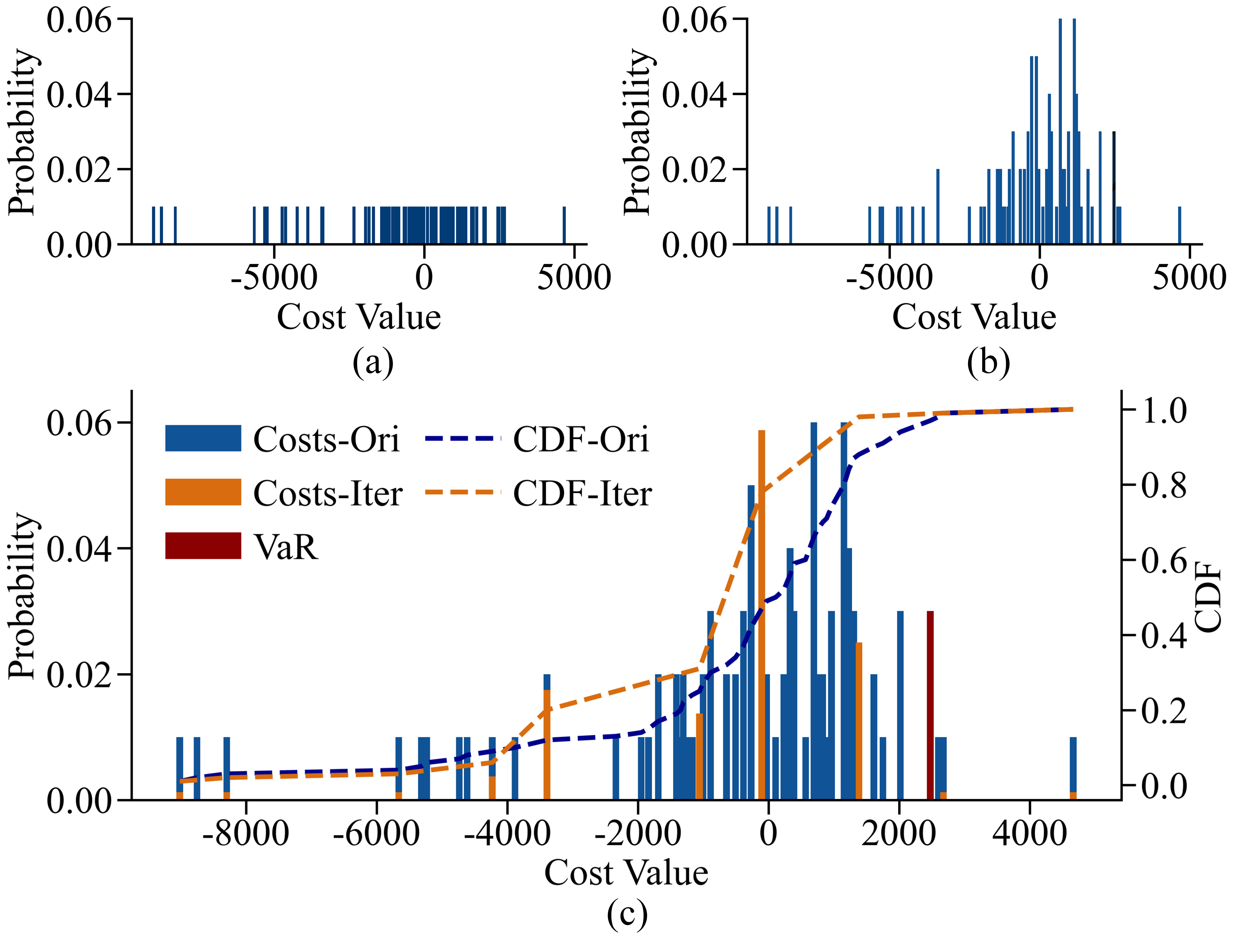}}
    \caption{Illustration example of the iterative scenario reduction process:
    (a) initial values of $\{f(z^\ast_{\bm{\zeta}},\xi_i)\}_{i=1}^N$;
    (b) aggregated values of $\{f(z^\ast_{\bm{\zeta}},\xi_i)\}_{i=1}^N$;
    (c) reduced values of $\bm{\zeta}$ after one iterative process.
    } 
    \label{fig:PDSR_iter_one}
\end{figure}

Fig.~\ref{fig:PDSR_iter_one} adopts a single iteration process to illustrate the key steps and effects of IPDSR.
Fig.~\ref{fig:PDSR_iter_one}(a) presents the initial distribution of the validation values $\{f(z^\ast_{\bm{\zeta}}, \xi_i)\}_{i=1}^N$. 
Fig.~\ref{fig:PDSR_iter_one}(b) shows the result after an effective objective aggregation, 
where the number of values is reduced from 100 to 50. 
This dimensionality reduction significantly alleviates the computational burden of the subsequent MIP problem~\eqref{eq:mip},
while preserving the essential characteristics of the original distribution.
Based on the aggregated data in Fig.~\ref{fig:PDSR_iter_one}(b),  
the MIP~\eqref{eq:mip} is then solved with approximation loss 0,
and the results are shown by the orange bars in Fig.~\ref{fig:PDSR_iter_one}(c).
Fig.~\ref{fig:PDSR_iter_one}(c) reveals that the cumulative distribution function (CDF) curve after the IPDSR iteration (orange line) 
closely aligns with the original distribution (blue line),
particularly in the right-tail region beyond the VaR threshold.
Notably, for clear illustration, the probabilities of Costs-iter are divided by 8.
This indicates that IPDSR effectively preserves critical risk information during scenario selection, 
maintaining tail structure and thereby enhancing the risk consistency and robustness of the resulting solution.

\begin{figure}[htbp] 
    \setlength{\abovecaptionskip}{-0.1cm}   
    \setlength{\belowcaptionskip}{-0.1cm}  
    \centerline{\includegraphics[width=0.95\columnwidth]{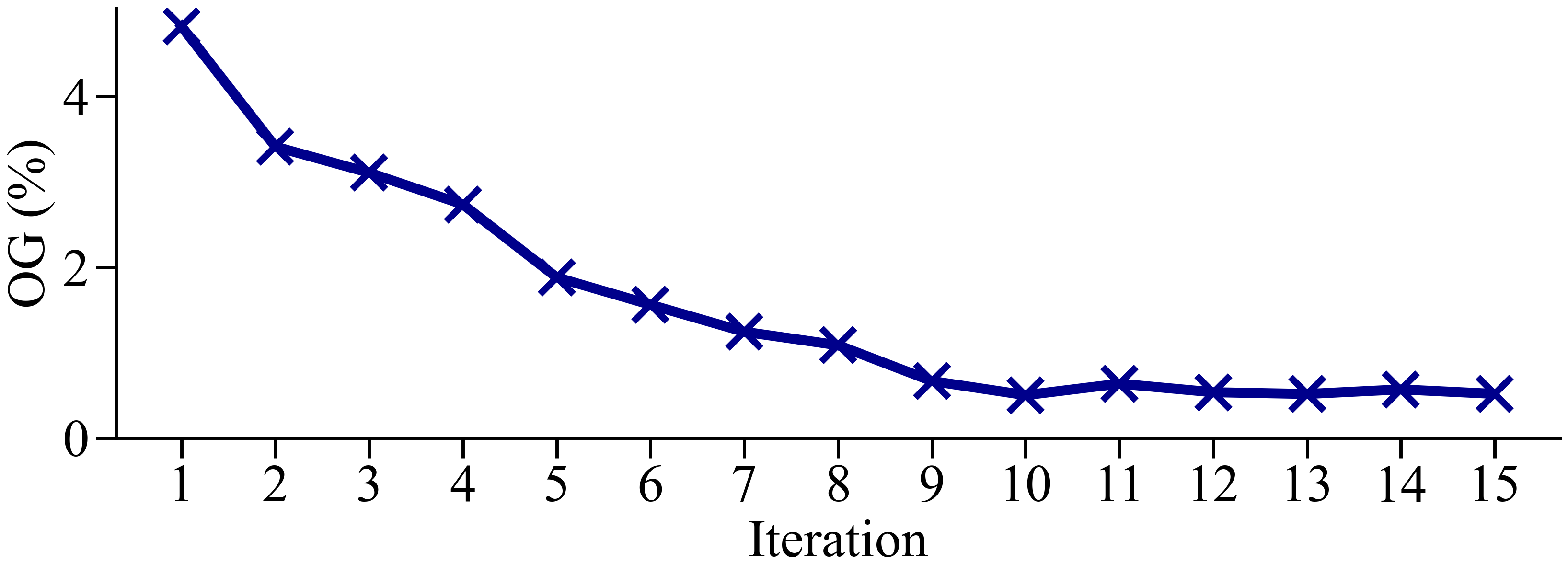}}
    \caption{The iterative curve of the OG value of IPDSR.}
    \label{fig:PDSR_iter_curve}
\end{figure}

Moreover,
we present the iterative curve of $OG(\%)$ in Fig.~\ref{fig:PDSR_iter_curve}.
The $OG(\%)$ value converges to a smaller value after a few iterations,
demonstrating the effectiveness of the iteration process in IPDSR.

Next, we utilize the evaluation indices of representative scenario effectiveness
to further validate the performance of IPDSR, as illustrated in Table~\ref{clustering eval tab}.
\begin{table}[htbp]
    \setlength{\abovecaptionskip}{-0.1cm}
    \setlength{\belowcaptionskip}{-0.1cm}
    \renewcommand\arraystretch{1.5} 
    \setlength{\tabcolsep}{0.1cm} 
    \caption{Evaluation results of representative scenario effectiveness\\ for $N=100$ and $K=10$}\label{clustering eval tab}
    \begin{center}
    \begin{tabular}{ccccccccccc}
    \toprule 
                               &$\zeta_1$ &$\zeta_2$ &$\zeta_3$ &$\zeta_4$ &$\zeta_5$ &$\zeta_6$&$\zeta_7$&$\zeta_8$&$\zeta_9$&$\zeta_{10}$ \\ 
    \midrule
    \textbf{$\omega_k$}         &$0.3 $&$0.18 $&$0.05 $&$0.05 $&$0.01 $&$0.06 $&$0.04 $&$0.16 $&$0.07 $&$0.08 $\\
    \textbf{$SE_{\zeta_k}(\%)$}  &$0.45 $&$0.49 $&$0.04 $&$0.02 $&$0.52 $&$0.09 $&$0.03 $&$0.35 $&$0.57 $&$0.04 $\\ 
    \bottomrule   
    \end{tabular}
    \end{center}
\end{table}	

Table~\ref{clustering eval tab} suggests that some scenarios of the representative scenario set has a relatively large value of $SE_{\zeta_k}(\%)$,
indicating their significant impact on the outcomes of the reduced problem. 
This result demonstrates the effectiveness of the proposed IPDSR framework in capturing the salient features 
of the original scenario set.

\subsection{Comparison with State-of-the-Art SR Methods}

To further demonstrate the benefits of IPDSR, 
we conduct a comparative analysis with $N=400$ and $K=20$.
DDSR methods including 
HC using Euclidean distance (HC-E),
HC using Wasserstein distance (HC-W), 
$K$-means using Euclidean distances (KM-E), 
$K$-means using dynamic time warping distance (KM-D), 
and Gaussian mixture model using Mahalanobis distance (G-M) are compared. 
Regarding other PDSR research, 
we also include the method developed in~\cite{SR_CVaR1}
focusing on iterative tail scenario reduction (ITSR) for comparison.
Recall that the PDSR in Section~\ref{subsec:PDSR} is not tractable for existing standard solvers, PDSR is not included for comparison.

The comparative indices include the SR performance indices and the computation efficiency indices.
We define the values fall beyond the VaR threshold of $\{f(z^\ast_{\bm{\xi}},\xi_i)\}_{i=1}^N$ as the ``worst-case'' scenarios.
The considered SR performance indices include:
the number of worst-case scenarios captured in the representative scenario set $(\varrho)$,
the $WD$ defined in~\eqref{eval:Wasserstein}, 
and the $OG(\%)$ defined in~\eqref{OG in eval index}. 
The considered computation efficiency indices include:
the time required to process the data input for clustering $(\tau_\mathrm{p})$,
and the time required to solve the clustering problem $(\tau_\mathrm{c})$.
The comparison results are presented in Table~\ref{tab: SRcompare}.
\begin{table}[htbp]
    \setlength{\abovecaptionskip}{-0.1cm}
    \setlength{\belowcaptionskip}{-0.1cm}
    \renewcommand\arraystretch{1.5} 
    \caption{Comparing SR methods for $N=400$ and $K=20$}\label{tab: SRcompare}
    \begin{center}
    \begin{tabular}{ccccccc}
    \toprule 
    Method & $\varrho$ & $WD$ & $OG(\%)$    &$\tau_\mathrm{p} (s)$ & $\tau_\mathrm{c} (s)$ \\ 
    \midrule
    Benchmark& $19$ & $0$  & $0$ & $-$ & $-$ \\
    IPDSR & $3$& $22.98$  & $0.73$ & $< 0.3$ & $<35$ \\
    PDSR& \multicolumn{5}{c}{Not tractable for existing standard solvers} \\
    ITSR~\cite{SR_CVaR1} & $9$& $162.87$  & $8.86$ & $< 0.3$ & $0.01$ \\
    G-M& $0$ & $48.3$  & $1.89$ & $-$ & $1.63$ \\
    KM-D& $0$ & $36.34$  & $2.74$ & $-$ & $12.7$ \\
    KM-E& $0$ & $42.95$  & $2.09$ & $-$ & $0.29$ \\
    HC-E& $0$ & $33.7$  & $1.37$ & $-$ & $0.02$ \\ 
    HC-W& $0$ & $90.77$  & $1.59$ & $8.8$ & $0.03$ \\ 
    \bottomrule  
    \end{tabular}
    \end{center}
\end{table}

\textit{\textbf{Optimality Gap:}} Table~\ref{tab: SRcompare} suggests that
IPDSR, with a small value of $OG(\%)=0.73\%$, 
significantly outperforms other SR methods. 
The DDSR methods, seeking for minimum statistical difference, 
fail to capture worst-case scenarios in their representative scenario sets.
Such limitations lead to a neglect of potential risks during the operation, 
and insufficient ability to make risk-averse decisions for potential future uncertainties, 
thus achieving relatively high $OG(\%)$.
This discrepancy highlights that the distribution similarity does not necessarily translate into approximated decision-making optimality in the problem space.
The ITSR method~\cite{SR_CVaR1} concentrates exclusively on tail scenarios and disregards the expected loss beyond the CVaR term, 
thereby incurring a significant OG bias and rendering it not suitable for the SBSO-CVaR formulation considered in this paper.
Compared to the above methods, 
IPDSR efficiently considers the potential impacts of the scenarios on the risk-averse problem, 
and includes three reasonable worst-case scenarios in the representative scenario set.
These comparative findings suggest that IPDSR exhibits superior accuracy in representing the optimality of original scenario set, 
thereby offering more reliable information for decision-making under uncertainty.


\textit{\textbf{Computational Efficiency:}} 
The Benchmark struggles with computational complexity when $N$ is large,
especially for complex SBSO-CVaR problems.
For example, for $N=400$ and $N=1000$ and a simple SBSO-CVaR problem as in Section~\ref{sec:formulation},
the Benchmark takes nearly 4 mins and 10 mins. 
Moreover, for larger $N$ and more complex SBSO-CVaR problems, the original problem becomes computationally intractable.
The DDSR methods, overcome computation bottlenecks, but at the price of potentially large $OG(\%)$ values.
IPDSR involve $N$ parallel computations of $\bm{F} = \{f(z^\ast_{\bm{\zeta}}, \xi_i)\}_{i=1}^N$ in each iteration
to obtain the validation values,
which can be efficiently executed in parallel.
In this paper, the computation time for each scenario-specific subproblem is $\tau_\text{s} < 0.3s$. 
For $N=400$ and $K=20$, after aggregating the $\bm{F}$ matrix into 100 values, 
the MIP clustering problem is solved consistently with MIP-gap with $\tau_\text{c} < 35s$. 
Thus, the IPDSR framework can effectively reduce computational complexity while achieving a small SR optimality gap.

\subsection{Scalability of IPDSR}
\begin{figure}[htbp] 
    \setlength{\abovecaptionskip}{-0.1cm}   
    \setlength{\belowcaptionskip}{-0.1cm}  
    \centerline{\includegraphics[width=0.95\columnwidth]{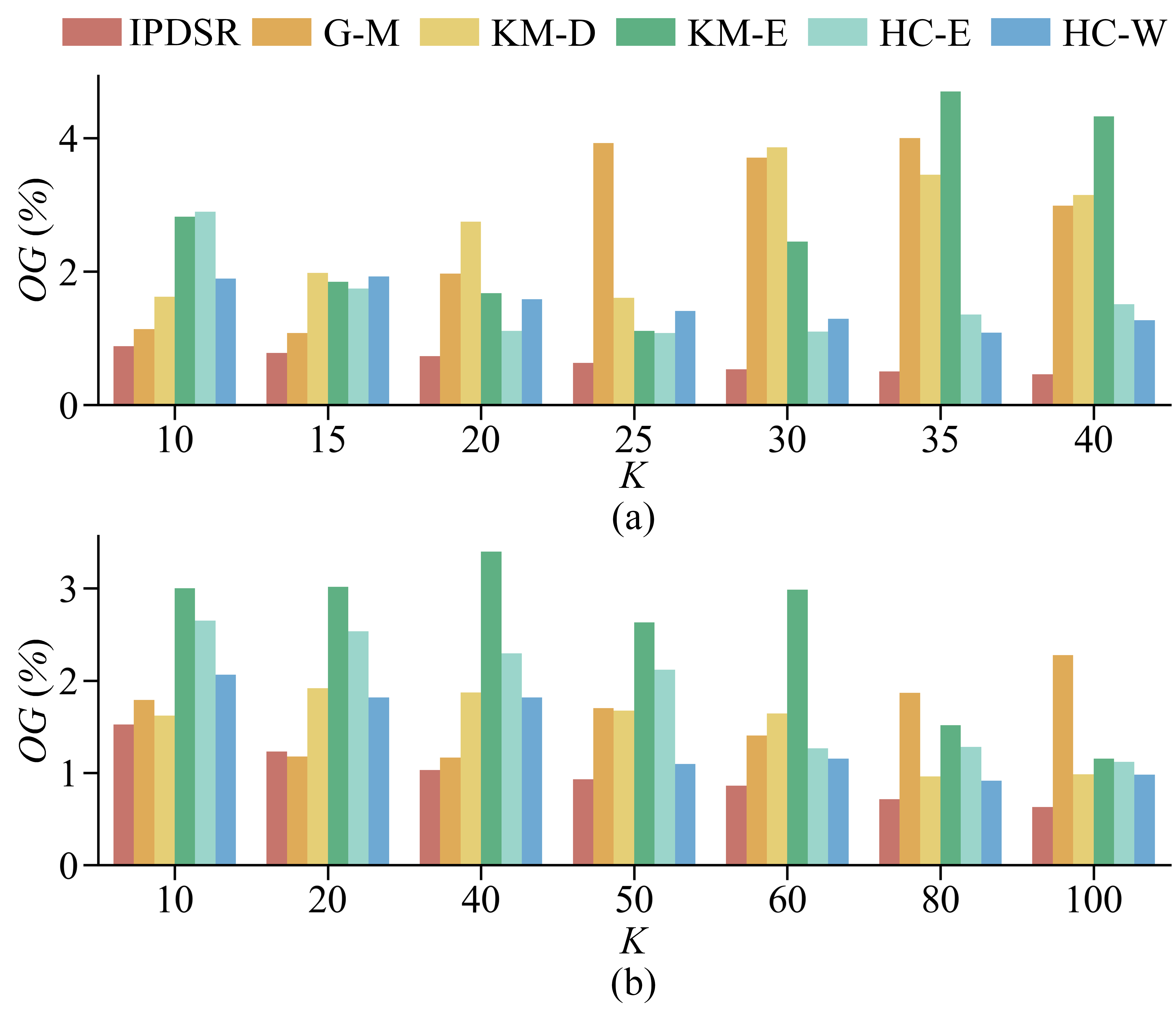}}
    \caption{Comparing results for different $N$ and $K$: (a) $N=400$, (b) $N=1000$.}
    \label{fig:SR_compare}
\end{figure}

We further compare the $OG(\%)$ results of multiple SR methods
under different $N$ and $K$, as illustrated in Fig.~\ref{fig:SR_compare}.
Note that the $OG(\%)$ of ITSR~\cite{SR_CVaR1} is all above 8\% and is thus not included.
In Fig.~\ref{fig:SR_compare},
we observe that IPDSR outperforms DDSR methods across all values of $N$ and $K$, 
as evidenced by its consistently lower $OG(\%)$. 
Besides,
as expected, increasing $K$ decreases the $OG(\%)$ for IPDSR,
as more representative scenarios can better capture the distribution features,
especially the right-tail scenarios. 
However, for DDSR methods, this inherent regular benefit is not present,
indicating that DDSR methods may capture some scenarios that are not impactful on the problem outcomes.
This problem-driven reduction mechanism ensures that the retained scenarios are truly representative with respect to decision-making relevance. 
Based on the representative scenarios from IPDSR, 
the SBSO-CVaR problem obtains superior operational strategies with high adaptability and effectiveness for future potential uncertainties. 
Moreover,
unlike SR for the risk-neutral SBSO problem~\cite{PDSR},
the SBSO-CVaR formulation requires a markedly larger set of representative scenarios 
to achieve a satisfactory $OG(\%)$. 
This is attributable to its pronounced emphasis on tail-risk, 
which can be adequately characterized only when the scenarios are sufficiently rich.

\section{Conclusion}\label{sec:conclusion}
In this paper, we propose an IPDSR framework for power system SBSO-CVaR problems, 
which fully incorporates the risk-averse problem structure into the SR process.
Specifically, we aim at minimizing the optimality gap,
and iteratively search for a better set of salient representative scenarios that can best approximate the 
validation distribution of the original scenario set and minimize the optimality gap. 
Moreover, we propose an MIP formulation to conduct the scenario partitioning and representative scenarios selection.
In the case study,
we visualize the nonlinear mapping relationship between the distribution space and the problem space.
By conducting SR in the problem space,
the presented IPDSR method obtains near-optimal approximation accuracy 
with just a few salient representative scenarios within acceptable computation time.
Moreover, a comprehensive comparative analysis with other SR methods is provided 
for different $N$ and $K$ values and demonstrates broadly the performance of our IPDSR method.

Future work will extend the method proposed in this paper to more types of risks like standard deviation, and further improve the efficiency.

\appendix\label{sec:appendix_PDSR}
In this appendix, we provide the detailed derivation process of the PDSR method
in Section~\ref{subsec:PDSR}.
including the mathematical derivation process of the upper bound of the OG,
the corresponding MIP formulation that minimizes the upper bound of the OG,
and the algorithmic steps.

\subsection{Derivation Process of the Upper Bound of OG}\label{subsec:appendix_PDSR}

For brief of mathematical notation, we denote  $G(z,\xi_i)$ and $G(z,\zeta_k)$ as follows:
\begin{subequations}
    \begin{align}
    &G(z,\xi_i) = f(z,\xi_i)  + \frac{\lambda}{1-\alpha} [ f(z,\xi_i) -v^{\bm{\xi}}_{\alpha} ]_{+}\\
    &G(z,\zeta_k) = f(z,\zeta_k)  + \frac{\lambda}{1-\alpha} [ f(z,\zeta_k) -v^{\bm{\zeta}}_{\alpha} ]_{+}
    \end{align}
\end{subequations}

Thus, we have $F(z,\bm{\xi}) = \sum_{i=1}^{N}\gamma_i G(z,\xi_i) + \lambda v^{\bm{\xi}}_{\alpha}$ 
and $F(z,\bm{\zeta}) = \sum_{k=1}^{K}\omega_k G(z,\zeta_k) + \lambda v^{\bm{\zeta}}_{\alpha}$.
Note that both the first and the last pair of~\eqref{eq:optimality gap} share a common expression of $F(z,\bm{\xi})-F(z,\bm{\zeta})$.
Since $\sum_{i=1}^{N}\gamma_i = \sum_{k=1}^{K} \sum_{i\in C_k}\gamma_i = \sum_{k=1}^{K}\omega_k$,
we have
\begin{align}
    &|F(z,\bm{\xi}) - F(z,\bm{\zeta})|\notag  \\
    =&\left|\sum_{k=1}^{K}\sum_{i\in C_k}\gamma_i G(z,\xi_i)+ \lambda v^{\bm{\xi}}_{\alpha} - \left(\sum_{k=1}^{K}\sum_{i\in C_k}\gamma_iG(z,\zeta_k)+ \lambda v^{\bm{\zeta}}_{\alpha}\right)\right|\notag \\
    =&\left| \sum_{k=1}^{K}\sum_{i\in C_k}\gamma_i\left(G(z,\xi_i)-G(z,\zeta_k) \right)+\lambda ( v^{\bm{\xi}}_{\alpha} - v^{\bm{\zeta}}_{\alpha})\right|\notag \\
    \leq&\sum_{k=1}^{K}\sum_{i\in C_k}\gamma_i| G(z,\xi_i)-G(z,\zeta_k)|+\lambda | v^{\bm{\xi}}_{\alpha} - v^{\bm{\zeta}}_{\alpha}|
\end{align}

Similar to Section~\ref{subsec:ISR},
We sort the scenarios in non-decreasing manner with respect to the validated objective function $f(z,\cdot)$.
Thus, we have
\begin{equation}
    \begin{split}
    &G(z,\xi_i) = f(z,\xi_i)  + \frac{\lambda}{1-\alpha} ( f(z,\xi_i) -f(z,\xi_{N^{\alpha}_{\xi}}) )\mathbf 1_{i>N^{\alpha}_{\bm{\xi}}}\\
    &G(z,\zeta_k) = f(z,\zeta_k)  + \frac{\lambda}{1-\alpha} ( f(z,\zeta_k) -f(z,\zeta_{K^{\alpha}_{\zeta}}) )\mathbf 1_{i>N^{\alpha}_{\bm{\zeta}}}
    \end{split}
\end{equation}
where $\mathbf 1_{\mathrm{condition}}$ is the indicator function, which equals 1 if the condition is true and 0 otherwise.

We give a uniform formulation of $|G(z,\xi_i)-G(z,\zeta_k)|$ as:
\begin{equation}
    \begin{split}
    &|G(z,\xi_i)-G(z,\zeta_k)|  \\
    =&|a(f(z,\xi_i)-f(z,\zeta_k)) + b(f(z,\xi_{N^{\alpha}_{\xi}})-f(z,\zeta_{K^{\alpha}_{\zeta}})) \\
    &+c(f(z,\xi_i)-f(z,\xi_{N^{\alpha}_{\xi}})) + d(f(z,\zeta_k)-f(z,\zeta_{K^{\alpha}_{\zeta}}))|\\
    \leq & |a| |f(z,\xi_i)-f(z,\zeta_k)| + |b| |f(z,\xi_{N^{\alpha}_{\xi}})-f(z,\zeta_{K^{\alpha}_{\zeta}})|\\
    &+|c| |f(z,\xi_i)-f(z,\xi_{N^{\alpha}_{\xi}})| + |d| |f(z,\zeta_k)-f(z,\zeta_{K^{\alpha}_{\zeta}})|\\
    \end{split}
\end{equation}
where $a$, $b$, $c$, and $d$ are coefficients summarized in Table~\ref{tab:cases}.
\begin{table}[!h]
\centering
\renewcommand{\arraystretch}{1.5} 
\caption{Values of $a$, $b$, $c$, and $d$ for different cases}
\begin{tabular}{ccccc}
\toprule
Case & $a$ & $b$ & $c$ & $d$ \\
\midrule
$i\leq N^{\alpha}_{\xi} \text{ and } k\leq K^{\alpha}_{\zeta}$ & 1 & 0 & 0 & 0 \\
$i> N^{\alpha}_{\xi} \text{ and } k\leq K^{\alpha}_{\zeta}$ & 1 & 0 & $\frac{\lambda}{1-\alpha}$ & 0 \\
$i\leq N^{\alpha}_{\xi} \text{ and } k> K^{\alpha}_{\zeta}$ & $1 + \frac{\lambda}{1-\alpha}$ & $-\frac{\lambda}{1-\alpha}$ & $-\frac{\lambda}{1-\alpha}$ & 0 \\
$i> N^{\alpha}_{\xi} \text{ and } k> K^{\alpha}_{\zeta}$ & $1 + \frac{\lambda}{1-\alpha}$ & $-\frac{\lambda}{1-\alpha}$ & 0 & 0 \\
\bottomrule
\end{tabular}
\label{tab:cases}
\end{table}

We introduce binary indicator variable $g_i$ to indicate $\mathbf 1_{i>N^{\alpha}_{\bm{\xi}}}$,
and $h_k$ to indicate $\mathbf 1_{k>K^{\alpha}_{\bm{\zeta}}}$.
Then we have
\begin{equation}
    \begin{split}
        &|G(z,\xi_i)-G(z,\zeta_k)| \\
        \leq &(h_k+(1 +\frac{\lambda}{1-\alpha})(1-h_k))|f(z,\xi_i)-f(z,\zeta_k)| \\
        + &\frac{\lambda}{1-\alpha}(1-h_k)|f(z,\xi_{N^{\alpha}_{\xi}})-f(z,\zeta_{K^{\alpha}_{\zeta}})|\\
        +&\frac{\lambda}{1-\alpha}|h_k - g_i||f(z,\xi_i)-f(z,\xi_{N^{\alpha}_{\xi}})|
    \end{split}
\end{equation}

Similar to our previous work~\cite{PDSR},
we define the problem-driven distance function $d(\xi_i,\zeta_k)$ as in~\eqref{eq:PDD}
and have the relation as in~\eqref{eq:lemma}.
Please refer to \cite{PDSR} for more details about~\eqref{eq:PDD} and~\eqref{eq:lemma}. 
\begin{align}
    &d(\xi_i,\zeta_k) = f(z_{\zeta_k}^\ast,\xi _i )-f(z_{\xi_i}^\ast,\xi_i)+f(z_{\xi_i}^\ast,\zeta _k )-f(z_{\zeta_k}^\ast,\zeta_k).\label{eq:PDD} \\
   & |f(z,\xi_i)-f(z,\zeta_k)|\leq h(\|z\Vert) d(\xi_i,\zeta_k)\text{,} \label{eq:lemma} 
\end{align}

Then, we have 
\begin{align} 
    &|F(z,\bm{\xi}) - F(z,\bm{\zeta})| \leq \sum_{k=1}^{K}\sum_{i\in C_k}\gamma_i|G(z,\xi_i)-G(z,\zeta_k)|+\lambda | v^{\bm{\xi}}_{\alpha} - v^{\bm{\zeta}}_{\alpha}|\notag\\
    &\leq h(\|z\Vert)(\sum_{k=1}^{K}\sum_{i\in C_k}\gamma_i( (h_k+(1 +\frac{\lambda}{1-\alpha})(1-h_k))d(\xi_i,\zeta_k) \notag\\
    &\ \ + (\frac{\lambda}{1-\alpha}(1-h_k)) d(\xi_{N^{\alpha}_{\bm{\xi}}},\zeta_{K^{\alpha}_{\bm{\zeta}}}) + \frac{\lambda}{1-\alpha}|h_k - g_i| d(\xi_i,\xi_{N^{\alpha}_{\bm{\xi}}}))\notag\\
    &\ \ +\lambda d(\xi_{N^{\alpha}_{\bm{\xi}}},\zeta_{K^{\alpha}_{\bm{\zeta}}}))
    \end{align} 

Since $Z$ is bounded, $\|z\Vert$ is well-defined, i.e, $\exists M\gg 1,~\|z\Vert \le M~\forall, z\in Z$.
Thus, we have the upper bound of OG as follows:
\begin{align}
        &OG\leq |F(z^\ast_{\bm{\zeta}},\bm{\xi}) - F(z^\ast_{\bm{\zeta}},\bm{\zeta})| + |F(z^\ast_{\bm{\xi}},\bm{\xi}) - F(z^\ast_{\bm{\xi}},\bm{\zeta})|\notag\\
        &\leq 2h(M)(\sum_{k=1}^{K}\sum_{i\in C_k}\gamma_i( (h_k+(1 +\frac{\lambda}{1-\alpha})(1-h_k))d(\xi_i,\zeta_k) \notag\\
    &\ \ + (\frac{\lambda}{1-\alpha}(1-h_k)) d(\xi_{N^{\alpha}_{\bm{\xi}}},\zeta_{K^{\alpha}_{\bm{\zeta}}}) + \frac{\lambda}{1-\alpha}|h_k - g_i| d(\xi_i,\xi_{N^{\alpha}_{\bm{\xi}}}))\notag\\
    &\ \ +\lambda d(\xi_{N^{\alpha}_{\bm{\xi}}},\zeta_{K^{\alpha}_{\bm{\zeta}}}))
\end{align}

\subsection{MIP Formulation}
Since $2h(M)$ is a constant, minimizing the following MIP formulation achieves the lowest upper bound of OG. 
\begin{subequations}\label{eq:full_model}
    \begin{align}
        \min
        &\ \sum_{j=1}^{N}\ l^{\mathrm{a}}_{j}+\sum_{j=1}^{N}\ l^{\mathrm{b}}_{j}+\sum_{j=1}^{N}l^{\mathrm{c}}_{j}+t \label{eq:obj_function} \\
        \text{s.t. } 
        & \sum_{i=1}^{N}v_{ij}\gamma_i(h_j+(1 +\frac{\lambda}{1-\alpha})(1-h_j))\tilde{F}_{ij} \le l^{\mathrm{a}}_{j}, \label{eq:constr_l1} \\
        & \sum_{i=1}^{N}v_{ij}\gamma_i(\frac{\lambda}{1-\alpha}(1-h_j))\sum_{p=1}^{N}\sum_{q=1}^{N}m_pn_q\tilde{F}_{pq} \le l^{\mathrm{b}}_{j},\label{eq:constr_l2} \\
        & \sum_{i=1}^{N}v_{ij}\gamma_i\frac{\lambda}{1-\alpha}|h_j - g_i|\sum_{p=1}^{N}m_p\tilde{F}_{ip} \le l^{\mathrm{c}}_{j}, \label{eq:constr_l3} \\
        & \lambda\sum_{p=1}^{N}\sum_{q=1}^{N}m_pn_q\tilde{F}_{pq} \leq t,\\
        & v_{ij} \leq u_{j},\ v_{jj} = u_{j},\ n_j \leq u_j,\label{eq:constr_selection1} \\
        & h_j \leq u_j,\ m_j \leq g_j,\ n_j \leq h_j, \label{eq:constr_selection2} \\
        & \sum\nolimits_{j=1}^{N}v_{ij} = 1, \label{eq:constr_assignment} \\
        & \sum\nolimits_{j=1}^{N}u_{j} = K, \label{eq:constr_K_clusters} \\
        & \sum\nolimits_{j=1}^{N}m_j = 1,\ \sum\nolimits_{j=1}^{N}n_j = 1, \label{eq:constr_anchor_centers} \\
        & \sum\nolimits_{i=1}^{N}g_i\gamma_i \leq \alpha, \label{var_g}\\
        & \sum\nolimits_{i=1}^{N}g_i\gamma_i + \gamma_k (1-g_k) \geq \alpha(1-g_k), \label{eq:constr_g_alpha} \\
        & \sum\nolimits_{j=1}^{N}h_j\sum_{i=1}^{N}v_{ij}\gamma_i \leq \alpha, \label{eq:constr_h_alpha_lower} \\
        & \sum_{j=1}^{N}h_j\sum_{i=1}^{N}v_{ij}\gamma_i + (u_k - h_k)\sum_{i=1}^{N}v_{ik}\gamma_i > \alpha (u_k - h_k), \label{eq:constr_h_alpha_upper} \\
        & \{v_{ij} ,\ u_{j} ,\ m_{j} ,\ n_{j} ,\ h_{k} ,\ g_{i}\} \in \{0,1\}. \label{eq:constr_binary_cont}
    \end{align}
\end{subequations}
where $\tilde{F}_{ij}:= f_{ij}-f_{jj} + f_{ji}-f_{ii}$ and $f_{ij}=f(z_{\xi_i}^\ast,\xi_j)$.
$m_j$ indicates whether scenario $\xi_j$ is the VaR scenario of $\bm{\xi}$,
and $n_j$ indicates whether scenario $\xi_j$ is the VaR scenario of $\bm{\zeta}$.
Constraint $\sum_{j=1}^{N}m_j=1$ ensures that there is only one VaR scenario in the original scenario set $\bm{\xi}$,
and $\sum_{j=1}^{N}n_j=1$ ensures that there is only one VaR scenario in the reduced scenario set $\bm{\zeta}$.
Constraint~\eqref{var_g}-\eqref{eq:constr_g_alpha} ensures that the selected VaR scenario is valid by satisfying (8).

However, the MIP formulation involves numerous binary variables and intricate constraints, 
making it computationally intractable for standard solvers.
We believe that this formulation can still provide valuable insights for SR in SBSO-CVaR problems,
and can be used as a theoretical basis for further research.

\subsection{Algorithm of PDSR}

The PDSR method is illustrated in \textbf{Algorithm~1}. 



\begin{algorithm}[htbp]\label{algo:PDSR}
\caption{Problem-Driven Scenario Reduction}
\SetAlgoLined
\SetEndCharOfAlgoLine{}
\KwIn{Scenario set $\boldsymbol{\xi}$ of $N$ scenarios and weights $\boldsymbol{\gamma}$.}
\KwOut{Scenario set $\boldsymbol{\zeta}$ of $K$ scenarios and weights $\boldsymbol{\omega}$.}
\SetKwBlock{StepOne}{Step 1 - Projection in Problem Space}{}
\SetKwBlock{StepTwo}{Step 2 - Clustering}{}
\SetKw{Parallel}{parallel}
\StepOne{
    Step 1.1: Initialize problem space matrix $\boldsymbol{F} = \mathbf{0}$.\;
    Step 1.2: Project the distribution space onto \newline problem space by:\;
    
    \For{$i=1$ \KwTo $N$ }{
        Solve the scenario-specific SO problem \newline $f(z,\xi_i )$ and obtain the optimal decision $z_{\xi_i}^\ast$. \;
           Set $f_{ii}=f(z_{\xi_i}^\ast,\xi_i)$.\;
        \For{$j=1$ \KwTo $N$, $i\neq j$ \Parallel}{
            Set $f_{ij}=f(z_{\xi_i}^\ast,\xi_j)$.
        }
    }
}
\StepTwo{
    Solve MIP in~\eqref{eq:full_model} to obtain the representative \newline scenario set $\boldsymbol{\zeta} $ with weights $\boldsymbol{\omega}$.\;
}
\end{algorithm}

\bibliographystyle{IEEEtran}
\bibliography{IEEEabrv,pscc2026}

\end{document}